\documentclass{amsart}
\usepackage{amstext,amssymb,amsthm,amsopn,newlfont,graphpap,graphics,graphicx,mathrsfs,enumitem}
\usepackage[parfill]{parskip}
\usepackage[noadjust]{cite}
\usepackage{epigraph}
\usepackage[colorlinks=true,
            linkcolor=red,
            urlcolor=blue,
            citecolor=magenta]{hyperref}
\usepackage{color}
\usepackage{mathrsfs}
\allowdisplaybreaks
\theoremstyle{plain}
\newtheorem{thm}{Theorem}[section]
\newtheorem*{thm*}{Theorem}

\newtheorem*{prop*}{Proposition}
\newtheorem{cor}{Corollary}[section]
\newtheorem*{cor*}{Corollary}

\newtheorem*{lem*}{Lemma}
\theoremstyle{definition}

\newtheorem*{defn*}{Definition}
\newtheorem{exmps}{Examples}[section]
\newtheorem*{exmps*}{Examples}
\newtheorem{exmp}[exmps]{Example}
\newtheorem*{exmp*}{Example}

\newtheorem*{exerc*}{Exercise}

\newtheorem{rems}{Remarks}[section]
\newtheorem*{rems*}{Remarks}
\newtheorem{rem}[rems]{Remark}
\newtheorem*{rem*}{Remark}
\newcommand{\N}{{\mathbb N}}
\newcommand{\Z}{{\mathbb Z}}
\newcommand{\R}{{\mathbb R}}
\newcommand{\C}{{\mathbb C}}

\newcommand{\emps}{\emptyset}

\renewcommand{\iff}{\: \Leftrightarrow\: }
\renewcommand{\bar}{\overline}

\numberwithin{equation}{section}

\DeclareMathOperator{\dist}{dist}

\DeclareMathOperator*{\esup}{\mbox{$E_A$}-ess\,sup}
\DeclareMathOperator*{\einf}{\mbox{$E_A$}-ess\,inf}

\begin{document}
\title[On spectral inclusion and mapping theorems]
{On spectral inclusion and mapping theorems\\
for scalar type spectral operators\\
and semigroups}
\author[Marat V. Markin]{Marat V. Markin}
\address{
Department of Mathematics\newline
California State University, Fresno\newline
5245 N. Backer Avenue, M/S PB 108\newline
Fresno, CA 93740-8001, USA
}
\email{mmarkin@csufresno.edu}
\dedicatory{To Michel L. Lapidus, a superb mathematician and an exceptional human being, in honor of his anniversary.}
\begin{abstract}
We establish spectral inclusion and mapping theorems for scalar type spectral operators, generalizing their counterparts for normal operators. Thereby, we extend a precise weak spectral mapping theorem, known to hold for $C_0$-semigroups of normal operators on complex Hilbert spaces, to the more general case of $C_0$-semigroups of scalar type spectral operators on complex Banach spaces. The finer spectrum structure is given itemized consideration.
\end{abstract}
\maketitle
\epigraph{\textit{I believe that mathematical reality lies outside us, that our function is to discover or observe it, and that the theorems which we prove, and which we describe grandiloquently as our ``creations," are simply the notes of our observations.}}{G.H. Hardy}

\section[Introduction]{Introduction}

We establish spectral inclusion and mapping theorems for scalar type spectral operators, generalizing their counterparts for normal operators (\cite[Theorem D.9]{HerLap} and \cite[Theorem D.11]{HerLap}). 

Thereby, we extend a precise \textit{weak spectral mapping theorem}, known to hold for $C_0$-semigroups of normal operators on complex Hilbert spaces, to the more general case of $C_0$-semigroups of scalar type spectral operators on complex Banach spaces, the a priory requirement of eventual norm continuity remaining superfluous (cf. \cite{arXiv:2103.05260}).

The finer spectrum structure is given itemized consideration.

\section[Preliminaries]{Preliminaries}

Here, we concisely outline essential preliminaries.

\subsection{Resolvent Set and Spectrum}\

For a closed linear operator $A$ in a complex Banach space $X$, the set
\[
\rho(A):=\left\{ \lambda\in \C \,\middle|\, \exists\, 
R(\lambda,A):={(A-\lambda I)}^{-1}\in L(X) \right\}
\]
($I$ is the \textit{identity operator} on $X$, $L(X)$ is the space of bounded linear operators on $X$) and its complement $\sigma(A):=\C\setminus \rho(A)$ are called its \textit{resolvent set} and \textit{spectrum}, respectively.

The spectrum is partitioned into pairwise disjoint subsets, called the \textit{point}, \textit{continuous}, and \textit{residual spectrum} of $A$, respectively, as follows:
\begin{equation*}
\begin{split}
& \sigma_p(A):=\left\{\lambda\in \C \,\middle|\,A-\lambda I\ \text{is \textit{not injective}, i.e., $\lambda$ is an \textit{eigenvalue} of $A$} \right\},\\
& \sigma_c(A):=\left\{\lambda\in \C \,\middle|\,A-\lambda I\ \text{is \textit{injective},
\textit{not surjective}, and $\overline{R(A-\lambda I)}=X$} \right\},\\
& \sigma_r(A):=\left\{\lambda\in \C \,\middle|\,A-\lambda I\ \text{is \textit{injective} and $\overline{R(A-\lambda I)}\neq X$} \right\}
\end{split}
\end{equation*}
($R(\cdot)$ is the \textit{range} of an operator and $\overline{\cdot}$ is the \textit{closure} of a set) (see, e.g., \cite{Dun-SchI,Markin2020EOT}).
 
\subsection{$C_0$-Semigroups}\

A $C_0$-semigroup $\left\{T(t) \right\}_{t\ge 0}$ on a complex Banach space $(X,\|\cdot\|)$ with generator $A$ is said to be subject to a \textit{weak spectral mapping theorem} if
\begin{equation*}\tag{WSMT}\label{WSMT}
\sigma(T(t))\setminus \{0\}=\overline{e^{t\sigma(A)}}\setminus \{0\},\ t\ge 0.
\end{equation*}

An \textit{eventually norm-continuous} $C_0$-semigroup $\left\{T(t) \right\}_{t\ge 0}$ on a complex Banach space, i.e., such that, for some $t_0>0$, the operator function
\[
[t_0,\infty)\ni t\mapsto T(t)\in L(X),
\]
is \textit{continuous} relative to the \textit{operator norm} $\|\cdot\|$ (here and henceforth, we use the same notation as for the norm on $X$), is subject to the following stronger version of \eqref{WSMT}:
\begin{equation*}\tag{SMT}\label{SMT}
\sigma(T(t))\setminus \{0\}=e^{t\sigma(A)},\ t\ge 0,
\end{equation*}
called a \textit{spectral mapping theorem} (see \cite[Proposition  V.$2.3$]{Engel-Nagel2006} and \cite[Theorem V.$2.8$]{Engel-Nagel2006}). The class of eventually norm-continuous $C_0$-semigroups encompasses $C_0$-semigroups with certain regularity properties, such as \textit{eventually compact} and \textit{eventually differentiable}, in particular \textit{analytic} and \textit{uniformly continuous} (see \cite[Section II.5]{Engel-Nagel2006}).

A $C_0$-semigroup $\left\{T(t) \right\}_{t\ge 0}$ (of \textit{normal operators}) on a complex Hilbert space generated by a \textit{normal operator} $A$ is subject to the following precise version of \textit{weak spectral mapping theorem} \eqref{WSMT}:
\begin{equation*}\tag{PWSMT}\label{PWSMT}
\sigma(T(t))=\overline{e^{t\sigma(A)}},\ t\ge 0,
\end{equation*}
\cite[Corollary V.$2.12$]{Engel-Nagel2006} without being a priori eventually norm-continuous. 

\subsection{Scalar Type Spectral Operators}\

A {\it scalar type spectral operator} is a densely defined closed linear operator $A$ in a complex Banach space with strongly $\sigma$-additive \textit{spectral measure} (the \textit{resolution of the identity}) $E_A(\cdot)$, which assigns to the Borel sets of the complex plane projection operators on $X$ and has the operator's \textit{spectrum} $\sigma(A)$ as its {\it support} \cite{Dunford1954,Survey58,Dun-SchIII}.

Associated with such an operator is the {\it Borel operational calculus}, assigning to each Borel measurable function $F:\sigma(A)\to \bar{\C}$ ($\bar{\C}:=\C\cup \{\infty\}$ is the extended complex plane) with $E_A\left(\left\{\lambda\in \C \,\middle|\, F(\lambda)=\infty\right\}\right)=0$ a scalar type spectral operator
\begin{equation*}
F(A):=\int\limits_{\sigma(A)} F(\lambda)\,dE_A(\lambda)
\end{equation*} 
in $X$, whose spectral measure is the image of $E_A(\cdot)$ under $F(\cdot)$, i.e.,
\begin{equation}\label{ism}
E_{F(A)}(\delta)=E_{A}(F^{-1}(\delta)),\ \delta\in \mathscr{B}(\C),
\end{equation}
($\mathscr{B}(\C)$ is the \textit{Borel $\sigma$-algebra} on $\C$), with
\begin{equation*}
A=\int\limits_{\sigma(A)} \lambda\,dE_A(\lambda)
\end{equation*}
\cite{Bade1954,Dunford1954,Survey58,Dun-SchIII}.

On a complex finite-dimensional Banach space, scalar type spectral operators are those linear operators, which furnish an \textit{eigenbasis} for the space, i.e., allow a diagonal matrix representation (see, e.g., \cite{Dunford1954,Survey58,Dun-SchIII}). 

In a complex Hilbert space, scalar type spectral operators are those that are similar to {\it normal operators} \cite{Wermer1954} (see also \cite{Lorch1939,Mackey1952}), the latter being the scalar type spectral operators for which the corresponding spectral measure projections are \textit{orthogonal} (see, e.g., \cite{Dun-SchII,Plesner}).

Various examples of scalar type spectral operators, including differential operators arising in the study of linear systems of partial differential equations, in particular perturbed Laplacians, can be found in \cite{Dun-SchIII}.

Due to its strong $\sigma$-additivity, the spectral measure is uniformly bounded, i.e.,
\begin{equation}\label{bounded}
\exists\, M\ge 1\ \forall\,\delta\in \mathscr{B}(\C):\ \|E_A(\delta)\|\le M
\end{equation}
(see, e.g., \cite{Dun-SchI}).


By \cite[Theorem XVIII.$2.11$ (c)]{Dun-SchIII}, for a Borel measurable function $F:\sigma(A)\to \bar{\C}$, the operator $F(A)$ is \textit{bounded} iff $F(\cdot)$ is $E_A$-\textit{essentially bounded}, i.e.,
\[
\esup_{\lambda\in \sigma(A)}|F(\lambda)|<\infty,
\]
in which case
\begin{equation}\label{boundedop}
\esup_{\lambda\in \sigma(A)}|F(\lambda)|\le \|F(A)\|
\le 4M\esup_{\lambda\in \sigma(A)}|F(\lambda)|,
\end{equation}
where $M\ge 1$ is from \eqref{bounded}.

For a scalar type spectral operator $A$, $\sigma(A)\neq \emps$,
\begin{equation}\label{pspg1}
\sigma_p(A)=\left\{ \lambda\in \C \,\middle|\,  E_A\left(\left\{ \lambda\right\}\right)\neq 	0\right\},
\end{equation}
with $E_A\left(\left\{ \lambda\right\}\right)X$ being the \textit{eigenspace} associated with
an eigenvalue $\lambda\in \sigma_p(A)$, i.e.,
\begin{equation}\label{pspg2}
E_A\left(\left\{ \lambda\right\}\right)X=\ker(A-\lambda I),\ \lambda \in \sigma_p(A),
\end{equation}
moreover,
\begin{equation}\label{cspg}
\sigma_c(A)=\left\{ \lambda\in \sigma(A) \,\middle|\,  
E_A\left(\left\{ \lambda\right\}\right)= 0\right\}
\end{equation}
and
\begin{equation}\label{rspg}
\sigma_r(A)=\emptyset
\end{equation}
\cite[Corollary 3.1]{Markin2017} (see also \cite{Markin2006}).

A scalar type spectral $C_0$-semigroup $\left\{T(t) \right\}_{t\ge 0}$ (i.e., a $C_0$-semigroup of scalar type spectral operators) on a complex Banach space $X$ is generated by a scalar type spectral operator \cite{Berkson1966,Panchapagesan1969}, which is the case \textit{iff}
\[
s(A)<\infty
\] 
with
\[
T(t)=e^{tA}:=\int\limits_{\sigma(A)} e^{t\lambda}\,dE_A(\lambda),\ t\ge 0,
\]
{\cite[Proposition $3.1$]{Markin2002(2)}}, the orbit maps of the semigroup
\[
T(t)f=e^{tA}f,\ t\ge 0,f\in X,
\] 
being the \textit{weak solutions} (also called the \textit{mild solutions}) of the associated abstract evolution equation
\begin{equation*}
y'(t)=Ay(t),\ t\ge 0,
\end{equation*}
\cite{Markin2002(1)} (see also \cite{Ball,Engel-Nagel2006}).

\section{Spectral Inclusion and Mapping Theorems}

\begin{thm}[Weak Spectral Inclusion and Mapping Theorem A.E.]\label{AEGSMT}\ \\
Let $A$ be a scalar type spectral operator $A$ in a complex Banach space $(X,\|\cdot\|)$ with spectral measure $E_A(\cdot)$.
\begin{enumerate}
\item If $F:\sigma(A)\to \C$ is a Borel measurable function, then
\begin{equation}\label{AEGWSMTincl}
\sigma(F(A))\subseteq \overline{F(\sigma(A)\setminus \sigma)}, 
\end{equation}
where $\sigma$ is an arbitrary Borel subset of $\sigma(A)$ for which $E_A(\sigma)=0$.
\item If furthermore the function $F:\sigma(A)\to \C$ is continuous on $\sigma(A)\setminus \sigma$, where $\sigma$ is a Borel subset of $\sigma(A)$ for which $E_A(\sigma)=0$, then
\begin{equation}\label{AEGWSMTeq}
\sigma(F(A))=\overline{F(\sigma(A)\setminus \sigma)}.
\end{equation}
\end{enumerate}
\end{thm} 

\begin{proof}
Let $F:\sigma(A)\to \C$ be a Borel measurable function, $\sigma$ be an arbitrary $E_A$-null subset of $\sigma(A)$, i.e., $E_A(\sigma)=0$, and
\[
\lambda\in \C\setminus \overline{F(\sigma(A)\setminus \sigma)}
\]
be arbitrary. Then
\[
\dist(\lambda,F(\sigma(A)\setminus \sigma)):=\inf_{\mu\in \sigma(A)\setminus \sigma}|F(\mu)-\lambda|>0
\]
(see, e.g., \cite{Markin2018EFA,Markin2020EOT}).

In view of $E_A(\sigma)=0$, by the properties of the Borel operational calculus (see \cite[Theorem XVIII.$2.11$ (e)]{Dun-SchIII}),
\begin{align*}
{(F(A)-\lambda I)}^{-1}
&=\int\limits_{\sigma(A)} \dfrac{1}{F(\mu)-\lambda}\,dE_A(\mu)
=\int\limits_{\sigma(A)} \dfrac{1}{F(\mu)-\lambda}
\chi_{\sigma(A)\setminus \sigma}(\mu)\,dE_A(\mu)
\\
&=:\int\limits_{\sigma(A)\setminus \sigma} \dfrac{1}{F(\mu)-\lambda}\,dE_A(\mu)
\end{align*}
($\chi_{\delta}(\cdot)$ is the \textit{characteristic function} of a set $\delta\subseteq \C$) is a \textit{bounded} linear operator on $X$, for which, by \eqref{boundedop}, when applied to the function 
\[
\sigma(A)\ni \mu \mapsto \dfrac{1}{F(\mu)-\lambda}
\chi_{\sigma(A)\setminus \sigma}(\mu)\in \C,
\]
we have:
\begin{align*}
\left\|{(F(A)-\lambda I)}^{-1}\right\|
&=\left\|\int\limits_{\sigma(A)} \dfrac{1}{F(\mu)-\lambda}
\chi_{\sigma(A)\setminus \sigma}(\mu)\,dE_A(\mu)\right\|
\\
&\le 4M\esup_{\mu\in \sigma(A)}\frac{1}{|F(\mu)-\lambda|}\chi_{\sigma(A)\setminus \sigma}(\mu)
\\
&=4M\esup_{\mu\in \sigma(A)\setminus \sigma}\frac{1}{|F(\mu)-\lambda|}
\\
&
\le 4M\sup_{\mu\in \sigma(A)\setminus \sigma}\dfrac{1}{|F(\mu)-\lambda|}
=4M\dfrac{1}{\inf_{\mu\in \sigma(A)\setminus \sigma}|F(\mu)-\lambda|}\\
&= \dfrac{4M}{\dist(\lambda,F(\sigma(A)\setminus \sigma))}<\infty,
\end{align*}
where $M\ge 1$ is from \eqref{bounded}.

Therefore, 
\[
\lambda \in \rho(F(A)).
\]

Thus, we have the inclusion
\begin{equation*}
\C\setminus \overline{F(\sigma(A)\setminus \sigma)}\subseteq \rho(F(A)),
\end{equation*}
or equivalently, 
\begin{equation*}
\sigma(F(A))\subseteq \overline{F(\sigma(A)\setminus \sigma)},
\end{equation*}
which completes the proof of part (1).

Now, suppose that the function $F:\sigma(A)\to \C$ is also continuous on $\sigma(A)\setminus \sigma$, where $\sigma$ is a Borel subset of $\sigma(A)$ for which $E_A(\sigma)=0$, and let
\[
\lambda\in \rho(F(A))
\]
be arbitrary.

By part (1), we have the inclusion
\begin{equation}\label{sincl1}
\sigma(F(A))\subseteq \overline{F(\sigma(A)\setminus \sigma)}.
\end{equation}

In view of $E_A(\sigma)=0$, by the properties of the Borel operational calculus (see \cite[Theorem XVIII.$2.11$ (e)]{Dun-SchIII}),
\begin{align*}
\int\limits_{\sigma(A)\setminus \sigma} \dfrac{1}{F(\mu)-\lambda}\,dE_A(\mu)
&:=\int\limits_{\sigma(A)} \dfrac{1}{F(\mu)-\lambda}
\chi_{\sigma(A)\setminus \sigma}(\mu)\,dE_A(\mu)
\\
&=\int\limits_{\sigma(A)} \dfrac{1}{F(\mu)-\lambda}\,dE_A(\mu)
\\
&={(F(A)-\lambda I)}^{-1}=:R(\lambda,F(A))\in L(X).
\end{align*}

Hence, by \eqref{boundedop}, when applied to the function 
\[
\sigma(A)\ni \mu \mapsto \dfrac{1}{F(\mu)-\lambda}
\chi_{\sigma(A)\setminus \sigma}(\mu)\in \C,
\]
we have:
\[
\esup_{\mu\in \sigma(A)\setminus \sigma}\frac{1}{|F(\mu)-\lambda|}
=\esup_{\mu\in \sigma(A)}\frac{1}{|F(\mu)-\lambda|}\chi_{\sigma(A)\setminus \sigma}(\mu)\le \|R(\lambda,A)\|<\infty,
\]
which, since
\[
\esup_{\mu\in \sigma(A)\setminus \sigma}\frac{1}{|F(\mu)-\lambda|}
=\frac{1}{\einf_{\mu\in \sigma(A)\setminus \sigma}|F(\mu)-\lambda|},
\]
implies that
\[
\einf_{\mu\in \sigma(A)\setminus \sigma}|F(\mu)-\lambda|>0.
\]

Therefore, there exists a Borel set $\delta\subseteq \sigma(A)\setminus \sigma$ such that
\begin{equation}\label{null1}
E_A(\delta)=0
\end{equation}
and
\[
\dist\left(\lambda,F\left(\left(\sigma(A)\setminus \sigma\right)\setminus \delta\right)\right):=\inf_{\mu\in \left(\sigma(A)\setminus \sigma\right)\setminus \delta}|F(\mu)-\lambda|>0,
\]
which implies that, for the open disk
\[
\Delta_r:=\left\{\mu\in \C\,\middle|\, |\mu-\lambda|<r\right\}
\]
centered at $\lambda$ with radius 
\[
r:=\dist\left(\lambda,F\left(\left(\sigma(A)\setminus \sigma\right)\setminus \delta\right)\right)>0,
\]
we have:
\begin{equation}\label{forall1}
\Delta_r\cap F\left(\left(\sigma(A)\setminus \sigma\right)\setminus \delta\right)= \emps.
\end{equation}

Assume that
\begin{equation}\label{exists}
\Delta_{r}\cap F(\delta)\neq \emps
\end{equation}
and let $\hat{F}(\cdot)$ be the restriction of $F(\cdot)$ to $\sigma(A)\setminus \sigma$. 

Jointly, \eqref{forall1} and \eqref{exists} imply that 
\begin{equation}\label{incl}
\emps \neq {\hat{F}}^{-1}\left(\Delta_{r}\right)\subseteq \delta.
\end{equation}

By the \textit{continuity} of $\hat{F}(\cdot)$ on $\sigma(A)\setminus \sigma$, the set
${\hat{F}}^{-1}\left(\Delta_{r}\right)$ is \textit{open} in $\sigma(A)\setminus \sigma$ (see, e.g., \cite{Markin2018EFA}), i.e.,
\begin{equation*}
\hat{F}^{-1}\left(\Delta_{r}\right)=\left(\sigma(A)\setminus \sigma\right)\cap \theta
=\left(\sigma(A)\cap \theta\right)\setminus \sigma,
\end{equation*}
where $\theta$ is an nonempty \textit{open} set in $\C$. 

Whence, by the properties of spectral measure (see, e.g., \cite{Survey58,Dun-SchIII}) and in view of $E_A(\sigma)=0$, we infer that
\begin{align*}
E_A\left(\hat{F}^{-1}\left(\Delta_{r}\right)\right)
&=E_A\left(\sigma(A)\cap \theta\right)
-E_A\left(\sigma(A)\cap \theta \cap \sigma \right)\\
&=E_A\left(\sigma(A)\cap \theta\right)
-E_A\left(\sigma(A)\cap \theta\right)E_A\left(\sigma \right)\\
&=E_A\left(\sigma(A)\cap \theta\right)
-E_A\left(\sigma(A)\cap \theta\right)0=E_A\left(\sigma(A)\cap \theta\right).
\end{align*}

The latter, since $\sigma(A)$ is the \textit{support} for the spectral measure $E_A(\cdot)$ (see Preliminaries) and $\theta$ is an nonempty open set in $\C$, implies that
\begin{equation}\label{nonzero}
E_A\left(\hat{F}^{-1}\left(\Delta_{r}\right)\right)=E_A\left(\sigma(A)\cap \theta\right)\neq 0.
\end{equation}

However, by the properties of spectral measure (see, e.g., \cite{Survey58,Dun-SchIII}) and in view of \eqref{incl} and \eqref{null1},
\begin{align*}
E_A\left(\hat{F}^{-1}\left(\Delta_{r}\right)\right)
&=E_A\left(\hat{F}^{-1}\left(\Delta_{r}\right)\cap \delta\right)\\
&=E_A\left(\hat{F}^{-1}\left(\Delta_{r}\right)\right)E_A\left(\delta\right)
=E_A\left(\hat{F}^{-1}\left(\Delta_{r}\right)\right)0=0.
\end{align*}

The obtained \textit{contradiction} with \eqref{nonzero} proves that assumption \eqref{exists} is false, and hence,
\begin{equation}\label{forall2}
\Delta_{r}\cap F(\delta)=\emps.
\end{equation}

Jointly, \eqref{forall1} and \eqref{forall2} imply
\begin{equation*}
\Delta_{r}\cap F(\sigma(A)\setminus \sigma)=\emps,
\end{equation*}
which leads to the conclusion that
\[
\lambda\in \C\setminus \overline{F(\sigma(A)\setminus \sigma)}.
\] 

Hence, we have the inclusion
\[
\rho(F(A))\subseteq \C\setminus \overline{F(\sigma(A)\setminus \sigma)},
\]
or equivalently,
\begin{equation}\label{sincl2}
\overline{F(\sigma(A)\setminus \sigma)}\subseteq \sigma(F(A)).
\end{equation}

Inclusions \eqref{sincl1} and \eqref{sincl2} jointly imply 
\begin{equation*}
\sigma(F(A))=\overline{F(\sigma(A)\setminus \sigma)},
\end{equation*}
which completes the proof of part (2), an thus, of the entire statement.
\end{proof} 

\begin{rems}\
\begin{itemize}
\item Theorem \ref{AEGSMT} generalizes \cite[Theorem D.11]{HerLap}, which covers the case a \textit{normal operator} $A$ in a complex Hilbert space and a Borel measurable function $F:\sigma(A)\to \C$ with a countable $E_A$-null set of discontinuities.
\item In \cite[Theorem $5.2$]{Bade1954}, equality \eqref{AEGWSMTeq} is stated for a \textit{spectral operator} $A$ (not necessarily of scalar type) in a complex Banach space (see, e.g., \cite{Dunford1954,Survey58,Dun-SchIII}) and a function $F:\sigma(A)\to \C$ analytic on $\sigma(A)$ with the exception of a finite subset $\sigma$ of poles for which $E_A(\sigma)=0$ and either analytic or having a pole at infinity (see also \cite[Theorem XVIII.$2.21$]{Dun-SchIII}).
\item As the following examples demonstrate, for inclusion \eqref{AEGWSMTincl} and equality \eqref{AEGWSMTeq} to hold, the condition $E_A(\sigma)=0$ is essential, and furthermore, without the requirement for the function
$F:\sigma(A)\to \C$ to be continuous on $\sigma(A)\setminus \sigma$, inclusion \eqref{AEGWSMTincl}, or even the inclusion
\[
\sigma(F(A))\subseteq F(\sigma(A)\setminus \sigma),
\]
may be strict.
\end{itemize}
\end{rems}

\begin{exmps}\label{exmpsAEGSMT}\
\begin{enumerate}
\item[1.] For the \textit{self-adjoint}, and hence, \textit{scalar type spectral} (see Preliminaries), operator 
\[
l_2\ni x:=\left(x_n\right)_{n\in \N}\mapsto Ax:=\left(0x_1,x_2,\frac{1}{2}x_3,\dots\right)\in l_2
\]
($\N:=\left\{ 1,2,\dots\right\}$ is the set of \textit{natural numbers})
of multiplication by the real-termed sequence 
\[
a_n:=
\begin{cases}
0,& n=1,\\
1/(n-1),&n\ge 2, 
\end{cases}
\]
on the complex Hilbert space $l_2$ of square-summable sequences with
\[
\sigma(A)=\sigma_p(A)=\left\{0\right\}\cup \left\{1/n\right\}_{n\in \N}
\]
(see, e.g., \cite{Markin2020EOT}), the Borel measurable function 
\[
F(\lambda):=
\begin{cases}
1/\lambda,&\lambda \neq 0,\\
0,& \lambda =0,
\end{cases}
\]
discontinuous at $0$, and the Borel subset $\sigma:=\left\{0\right\}\subseteq \sigma(A)$ with $E_A\left(\sigma\right)\neq 0$ (see \eqref{pspg1}), $F(A)$ is the operator of multiplication by the sequence 
\[
F(a_n):=
\begin{cases}
0,& n=1,\\
n-1,&n\ge 2, 
\end{cases}
\]
i.e.,
\[
l_2\supseteq D\left(F(A)\right)\ni x:=\left(x_n\right)_{n\in \N}\mapsto F(A)x:=\left(0x_1,x_2,2x_3,\dots\right)\in l_2,
\]
where
\[
D\left(F(A)\right)=\left\{ \left(x_n\right)_{n\in \N}\in l_2\,\middle|\, \left(0x_1,x_2,2x_3,\dots\right)\in l_2 \right\}
\]
($D(\cdot)$ is the \textit{domain} of an operator), with 
\[
\sigma(F(A))=\sigma_p(F(A))=\Z_+
\]
($\Z_+:=\left\{0,1,2,\dots\right\}$ is the set of \textit{nonnegative integers}).

However,
\[
\sigma\left(F(A)\right)\not\subseteq \overline{F\left(\sigma(A)\setminus \sigma\right)}=F\left(\sigma(A)\setminus \left\{0\right\}\right)=\N.
\]
\item[2.] For the \textit{self-adjoint} operator
\[
L_2(\R)\supseteq D(A)\ni f\mapsto [Af](x):=xf(x)\in L_2(\R)
\]
of multiplication by the independent variable in the complex Hilbert space $L_2(\R)$, where
\[
D(A):=\left\{f\in L_2(\R)\,\middle|\, \int_0^\infty|xf(x)|^2\,dx<\infty  \right\}
\]
($f(\cdot)$ is a representative of an equivalence class $f\in L_2(\R)$), with
\[
\sigma(A)=\sigma_c(A)=\R
\]
(see, e.g., \cite{Akh-Glaz}), the Borel measurable function 
$\chi_{\left\{ 0\right\}}(\cdot)$ discontinuous at $0$, and the Borel subset $\sigma:=\emps\subseteq \sigma(A)$, 
\[
\chi_{\left\{ 0\right\}}(A)=E_A(\left\{ 0\right\})=0
\]
(see \eqref{cspg}) (see, e.g., \cite{Dun-SchII,Plesner}) with 
\[
\sigma\left(\chi_{\left\{ 0\right\}}(A)\right)= \sigma_p\left(\chi_{\left\{ 0\right\}}(A)\right)
=\left\{ 0\right\}.
\]

However,
\[
\sigma\left(\chi_{\left\{ 0\right\}}(A)\right)\neq \overline{\chi_{\left\{ 0\right\}}\left(\sigma(A)\right)}=\chi_{\left\{ 0\right\}}\left(\sigma(A)\right)=\left\{0,1 \right\}.
\]
\end{enumerate}
\end{exmps}

\begin{rem}\label{remexmpsAEGSMT}
As is seen from the prior example, for any $\lambda\in \R\setminus \left\{0\right\}$, $\lambda\in \sigma_c(A)$, but $\chi_{\left\{ 0\right\}}(\lambda)=0\in \sigma_p\left(\chi_{\left\{ 0\right\}}(A)\right)$ and $0\in \sigma_c(A)$, but $\chi_{\left\{ 0\right\}}(0)=1\in \rho\left(\chi_{\left\{ 0\right\}}(A)\right)$.
\end{rem}

The subsequent statement, being the particular case of Theorem \ref{AEGSMT} for $\sigma:=\emps$, generalizes \cite[Theorem D.9]{HerLap}, its counterpart for normal operators (cf. also \cite[Lemma 4, Theorem 2]{Ricker1985}).

\begin{thm}[Weak Spectral Inclusion and Mapping Theorem]\label{GSMT}\ \\
Let $A$ be a scalar type spectral operator $A$ in a complex Banach space.
\begin{enumerate}
\item If $F:\sigma(A)\to \C$ is a Borel measurable function, then
\begin{equation}\label{GWSMTincl}
\sigma(F(A))\subseteq \overline{F(\sigma(A))}.
\end{equation}
\item If $F:\sigma(A)\to \C$ is a continuous function, then
\begin{equation}\label{GWSMTeq}
\sigma(F(A))=\overline{F(\sigma(A))}.
\end{equation}
\end{enumerate}
\end{thm} 

\begin{rem}
As is seen from Examples \ref{exmpsAEGSMT}, without the requirement of continuity for the Borel measurable function $F:\sigma(A)\to \C$, inclusion \eqref{GWSMTincl} may be strict even when the set of discontinuities is an $E_A$-null set. 
\end{rem}

\section{Finer Spectrum Structure}

The finer spectrum structure is governed by the following statement.

\begin{thm}[Finer Spectrum Structure]\label{FST}\ \\
Let $A$ be a scalar type spectral operator in a complex Banach space with spectral measure $E_A(\cdot)$.
\begin{enumerate}
\item  If $F:\sigma(A)\to \C$ is a Borel measurable function, then
\begin{equation}\label{pspA}
\sigma_p(F(A)) \supseteq F(\sigma_p(A)) \quad \hfill \text{(Spectral Inclusion Theorem for Point Spectrum)}.
\end{equation}
\item If $F:\sigma(A)\to \C$ is a Borel measurable function which is also injective, then
\begin{equation}\label{pspAA}
\sigma_p(F(A))=F(\sigma_p(A))
\quad \hfill \text{(Spectral Mapping Theorem for Point Spectrum)},
\end{equation}
with the operators $A$ and $F(A)$ sharing the corresponding spectral measure projections, i.e.,
\begin{equation}\label{smp}
E_{F(A)}\left(\left\{F(\lambda)\right\}\right)=E_A(\left\{ \lambda\right\}),\ \lambda\in \sigma_p(A),
\end{equation}
and hence,  the corresponding eigenspaces, i.e.,
\begin{equation}\label{ees}
\ker\left(F(A)-F(\lambda)I\right)=\ker\left(A -\lambda I\right),\ \lambda\in \sigma_p(A),
\end{equation}
and 
\begin{equation}\label{cspAA}
\rho(F(A))\cup \sigma_c(F(A)) \supseteq F(\sigma_c(A)).
\end{equation}

In particular, if $\sigma_p(A)=\emptyset$, then $\sigma_p(F(A))=\emps$.
\item If furthermore $F:\sigma(A)\to \C$ is a continuous and injective function,
\begin{equation}\label{cspAAA}
\sigma_c(F(A)) = F(\sigma_c(A))\cup \left(\overline{F(\sigma(A))}\setminus F(\sigma(A))\right).
\end{equation}

In particular, if the set $F(\sigma(A))$ is closed, as is the case when the operator $A$ is bounded, and hence, so is $F(A)$,
\begin{equation}\label{cspAAAA}
\sigma_c(F(A)) = F(\sigma_c(A))
\ \text{(Spectral Mapping Thm for Continuous Spectrum)}.
\end{equation}
\end{enumerate}
\end{thm} 

\begin{proof}
Let $\lambda \in \sigma_p(A)$ be arbitrary. Then
\[
E_A(\left\{ \lambda\right\})\neq 0
\]
(see \eqref{pspg1}).

In view of the fact that the spectral measure $E_{F(A)}(\cdot)$ of $F(A)$ is defined by \eqref{ism} (see \cite[Theorem $3.3$]{Bade1954}) and by the properties of spectral measure (see, e.g., \cite{Survey58,Dun-SchIII}),
\begin{align*}
E_A(\left\{ \lambda\right\})E_{F(A)}\left(\left\{F(\lambda)\right\}\right)
&=E_A(\left\{ \lambda\right\})E_A\left(F^{-1}\left(\left\{F(\lambda)\right\}\right)\right)\\
&=E_A\left(\left\{ \lambda\right\}\cap F^{-1}\left(\left\{F(\lambda)\right\}\right)\right)
=E_A(\left\{ \lambda\right\})\neq 0.
\end{align*}

Whence, we conclude that
\[
E_{F(A)}\left(\left\{F(\lambda)\right\}\right)\neq 0,
\]
and hence,
\[
F(\lambda)\in \sigma_p(F(A)) 
\]
(see \eqref{pspg1}). 

Thus, inclusion \eqref{pspA} holds and the proof of part (1) is complete.

Now, suppose that the function $F:\sigma(A)\to \C$ is \textit{injective}.  Then there exists an \textit{inverse} $F^{-1}(\cdot)$.

By \eqref{pspg1} and \eqref{ism},
\begin{equation}\label{pspeq}
\begin{aligned}
\mu \in \sigma_p\left(F(A)\right) &\iff E_{F(A)}\left(\left\{\mu \right\}\right)\neq 0
\iff E_A\left(F^{-1}\left(\left\{\mu \right\}\right)\right)\neq 0 \\
&\iff F^{-1}\left(\left\{\mu \right\}\right)\neq \emps\ \text{and}\ \lambda:=F^{-1}(\mu)\in \sigma_p(A),
\end{aligned}
\end{equation}
which, by part (1), proves equality \eqref{pspAA}, with equalities \eqref{smp} and \eqref{ees} following by \eqref{ism} and \eqref{pspg2}.

In view of the \textit{injectivity} of the function $F:\sigma(A)\to \C$, equality \eqref{pspAA}
and the fact that
\begin{equation}\label{rspe}
\sigma_r(F(A))=\emps
\end{equation}
imply inclusion \eqref{cspAA}.

This completes the proof of part (2).

Finally, suppose that the function $F:\sigma(A)\to \C$ is \textit{continuous} and \textit{injective}. 

By the \textit{Weak Spectral Inclusion and Mapping Theorem} (Theorem \ref{GSMT}), the \textit{continuity} of the function $F:\sigma(A)\to \C$ implies that
\begin{equation}\label{wsmt}
\sigma(F(A))=\overline{F(\sigma(A))},
\end{equation}
or equivalently,
\[
\rho(F(A))=\C\setminus\overline{F(\sigma(A))}.
\]

By part (2), the \textit{injectivity} of the function $F:\sigma(A)\to \C$ implies equality \eqref{pspAA} and inclusion \eqref{cspAA}.

In view of equality \eqref{wsmt}, inclusion \eqref{cspAA} turns into the inclusion
\begin{equation}\label{cspincl1}
\sigma_c(F(A)) \supseteq F(\sigma_c(A)).
\end{equation}

Further, by equalities \eqref{wsmt}, \eqref{pspAA}, and \eqref{rspe}, we infer that
\begin{equation}\label{cspincl2}
\sigma_c(F(A)) \supseteq \overline{F(\sigma(A))}\setminus F(\sigma(A)).
\end{equation}

Thus, inclusions \eqref{cspincl1} and \eqref{cspincl2} imply the inclusion
\begin{equation}\label{cspincl3}
\sigma_c(F(A)) \supseteq F(\sigma_c(A))\cup \left(\overline{F(\sigma(A))}\setminus F(\sigma(A))\right).
\end{equation}

On the other hand, in view of equalities \eqref{wsmt}, \eqref{pspAA}, and \eqref{rspe}, we also have the converse inclusion
\begin{equation}\label{cspincl4}
\sigma_c(F(A)) \subseteq F(\sigma_c(A))\cup \left(\overline{F(\sigma(A))}\setminus F(\sigma(A))\right).
\end{equation}

Inclusions \eqref{cspincl3} and \eqref{cspincl4} jointly imply equality \eqref{cspAAA}.

This completes the proof of part (3), and thus, of the entire statement.
\end{proof}

\begin{samepage}
\begin{rems}\
\begin{itemize}
\item As is seen from Examples \ref{exmpsAEGSMT} (see Remark \ref{remexmpsAEGSMT}), without the requirement of injectivity for the Borel measurable function $F:\sigma(A)\to \C$, inclusion \eqref{pspA} may be strict and inclusion \eqref{cspAA} need not hold. 
\item The following example shows that, with the requirement of injectivity but without the requirement of continuity for the Borel measurable function $F:\sigma(A)\to \C$, equality \eqref{cspAAA} need not hold. 
\end{itemize}
\end{rems}
\end{samepage}

\begin{exmp}\label{exmpFST}
For the \textit{self-adjoint} operator $A$ of multiplication by the independent variable in the complex Hilbert space $L_2(\R)$ from Examples \ref{exmpsAEGSMT} and the real-valued Borel measurable function
\[
F(\lambda):=
\begin{cases}
e^\lambda,&\lambda \neq 0,\\
-1,& \lambda =0,
\end{cases}
\] 
on $\sigma(A)=\sigma_c(A)=\R$, which is injective but discontinuous at $0$ with 
$E_A(\left\{ 0\right\})=0$ (see \eqref{cspg}), $F(A)$ is the self-adjoint operator of multiplication by $F(\cdot)$.

By Theorem \ref{FST} (part (2)), which implies that $\sigma_p(F(A))=\emps$, and Theorem \ref{AEGSMT},
\[
\sigma_c(F(A))=\sigma(F(A))=\overline{F(\sigma(A)\setminus \left\{ 0\right\})}
=\overline{F(\R\setminus \left\{ 0\right\})}=[0,\infty).
\]

Thus, for any $\lambda\in \R\setminus \left\{0\right\}$, $\lambda\in \sigma_c(A)$ and $F(\lambda)=e^\lambda\in \sigma_c\left(F(A)\right)$ and, for $0\in \sigma_c(A)$, $F(0)=-1\in \rho\left(F(A)\right)$.
\end{exmp}

\begin{rem}
As follows from Examples \ref{exmpsAEGSMT} and Example \ref{exmpFST}, 
generally, for a scalar type spectral operator $A$ in a complex Banach space and a Borel measurable function $F:\sigma(A)\to \C$, 
\begin{equation*}
\rho(F(A))\cup \sigma_c(F(A))\cup \sigma_p(F(A)) \supseteq F(\sigma_c(A)).
\end{equation*}
\end{rem}

\section{Scalar Type Spectral $C_0$-Semigroups}

Considering that, for a $C_0$-semigroup $\left\{T(t) \right\}_{t\ge 0}$ (of scalar type spectral operators) on a complex Banach space generated by a scalar type spectral operator $A$ with spectral measure $E_A(\cdot)$,
\[
T(t)=e^{tA}:=\int\limits_{\sigma(A)} e^{t\lambda}\,dE_A(\lambda),\ t\ge 0,
\]
(see Preliminaries) and applying for each $t\ge 0$ the \textit{Weak Spectral Inclusion and Mapping Theorem} (Theorem \ref{GSMT}) relative to the continuous exponential function 
\[
F_t(\lambda):=e^{t\lambda},\ \lambda\in \sigma(A),
\]
we obtain the following generalization of \textit{precise weak spectral mapping theorem} \eqref{PWSMT}, known to hold for $C_0$-semigroups of normal operators (see Prelimina\-ries).

\begin{thm}[Precise Weak Spectral Mapping Theorem]\label{PWSMThm}\ \\
A $C_0$-semigroup $\left\{T(t) \right\}_{t\ge 0}$ (of scalar type spectral operators) on a complex Banach space generated by a scalar type spectral operator $A$ is subject to precise weak spectral mapping theorem \eqref{PWSMT}.
\end{thm} 

The finer spectrum structure for scalar type spectral $C_0$-semigroups is governed by the following statement.

\begin{thm}[Finer Spectrum Structure]\ \\
Let  $\left\{T(t) \right\}_{t\ge 0}$  be a $C_0$-semigroup (of scalar type spectral operators) on a complex Banach space generated by a scalar type spectral operator $A$. Then
\begin{align}
\sigma_p(T(t)) &=e^{t\sigma_p(A)},\, t\ge 0,\, \text{(Spectral Mapping Theorem for Point Spectrum)},\label{psp}\\ 
\sigma_r(T(t))&=\emptyset,\ t\ge 0,\label{rsp}\\
\sigma_c(T(t))&=\overline{e^{t\sigma(A)}}\setminus e^{t\sigma_p(A)},\ t\ge 0.\label{csp}
\end{align}

If furthermore, for some $t>0$, the restriction of the exponential function $e^{t\cdot}$ to $\sigma(A)$ is injective, then
\begin{equation}\label{cspi}
\sigma_c(T(t)) =e^{t\sigma_c(A)}\cup \left(\overline{e^{t\sigma(A)}}\setminus e^{t\sigma(A)}\right).
\end{equation}

In particular, if the set $e^{t\sigma(A)}$ is closed, as is the case when the operator $A$ is bounded,
\begin{equation}\label{cspii}
\sigma_c(T(t)) =e^{t\sigma_c(A)}
\ \text{(Spectral Mapping Theorem for Continuous Spectrum)}.
\end{equation}
\end{thm} 

\begin{proof}
By \cite[Theorem V.$2.6$]{Engel-Nagel2006}, 
\begin{equation}\label{pspsg}
\sigma_p(T(t))\setminus \left\{0\right\} =e^{t\sigma_p(A)},\ t\ge 0,
\end{equation}
even without the assumption of scalar type spectrality for the generator.

Since, in view of
\[
T(t)=e^{tA},\ t\ge 0,
\]
(see Preliminaries), by the properties of the Borel operational calsulus (see \cite[Theorem XVIII.$2.11$ (h)]{Dun-SchIII}), for any $t\ge 0$, there exists the inverse
\[
{T(t)}^{-1}={\left(e^{tA}\right)}^{-1}=e^{-tA},
\]
we infer that
\begin{equation}\label{0nev}
0\notin \sigma_p(T(t)),\ t\ge 0.
\end{equation}

From \eqref{pspsg} and \eqref{0nev}, we infer that equality \eqref{psp} holds. 

Equality \eqref{rsp} holds by \eqref{rspg}.

Equality \eqref{csp} follows by Theorem \ref{PWSMThm} from equalities \eqref{psp} and \eqref{rsp}.

If furthermore, for some $t>0$, the restriction of the exponential function $e^{t\cdot}$ to $\sigma(A)$ is injective, then, by Theorem \ref{FST} (part (3)), we arrive at equality \eqref{cspi}.
\end{proof} 

\section{The Case of Normal Operators and $C_0$-Semigroups}

For the case of normal operators and $C_0$-semigroups (see Preliminaries), we arrive at the following corollaries of the corresponding statements.

\begin{cor}[Weak Spectral Inclusion and Mapping Theorem A.E.]\ \\
Let $A$ be a normal operator $A$ in a complex Hilbert space with spectral measure $E_A(\cdot)$.
\begin{enumerate}
\item If $F:\sigma(A)\to \C$ is a Borel measurable function, then
\begin{equation*}
\sigma(F(A))\subseteq \overline{F(\sigma(A)\setminus \sigma)}, 
\end{equation*}
where $\sigma$ is an arbitrary Borel subset of $\sigma(A)$ for which $E_A(\sigma)=0$.
\item If furthermore the function $F:\sigma(A)\to \C$ is continuous on $\sigma(A)\setminus \sigma$, where $\sigma$ is a Borel subset of $\sigma(A)$ for which $E_A(\sigma)=0$, then
\begin{equation*}
\sigma(F(A))=\overline{F(\sigma(A)\setminus \sigma)}.
\end{equation*}
\end{enumerate}
\end{cor} 

Cf. \cite[Theorem D.11]{HerLap} (see also \cite[Theorem D.12]{HerLap}).

\begin{cor}[Weak Spectral Inclusion and Mapping Theorem]\ \\
Let $A$ be a normal operator $A$ in a complex Hilbert space.
\begin{enumerate}
\item If $F:\sigma(A)\to \C$ is a Borel measurable function, then
\begin{equation*}
\sigma(F(A))\subseteq \overline{F(\sigma(A))}.
\end{equation*}
\item If $F:\sigma(A)\to \C$ is a continuous function, then
\begin{equation*}
\sigma(F(A))=\overline{F(\sigma(A))}.
\end{equation*}
\end{enumerate}
\end{cor} 

Cf. \cite[Theorem D.9]{HerLap}.

\begin{cor}[Finer Spectrum Structure]\ \\
Let $A$ be a normal operator $A$ in a complex Hilbert space with spectral measure $E_A(\cdot)$.
\begin{enumerate}
\item  If $F:\sigma(A)\to \C$ is a Borel measurable function, then
\begin{equation*}
\sigma_p(F(A)) \supseteq F(\sigma_p(A)) \quad \hfill \text{(Spectral Inclusion Theorem for Point Spectrum)}.
\end{equation*}
\item If $F:\sigma(A)\to \C$ is a Borel measurable function which is also injective, then
\begin{equation*}
\sigma_p(F(A))=F(\sigma_p(A))
\quad \hfill \text{(Spectral Mapping Theorem for Point Spectrum)},
\end{equation*}
with the operators $A$ and $F(A)$ sharing the corresponding spectral measure projections, i.e.,
\begin{equation*}
E_{F(A)}\left(\left\{F(\lambda)\right\}\right)=E_A(\left\{ \lambda\right\}),\ \lambda\in \sigma_p(A),
\end{equation*}
and hence,  the corresponding eigenspaces, i.e.,
\begin{equation*}
\ker\left(F(A)-F(\lambda)I\right)=\ker\left(A -\lambda I\right),\ \lambda\in \sigma_p(A),
\end{equation*}
and 
\begin{equation*}
\rho(F(A))\cup \sigma_c(F(A)) \supseteq F(\sigma_c(A)).
\end{equation*}

In particular, if $\sigma_p(A)=\emptyset$, then $\sigma_p(F(A))=\emps$.
\item If furthermore $F:\sigma(A)\to \C$ is a continuous and injective function,
\begin{equation*}
\sigma_c(F(A)) = F(\sigma_c(A))\cup \left(\overline{F(\sigma(A))}\setminus F(\sigma(A))\right).
\end{equation*}

In particular, if the set $F(\sigma(A))$ is closed, as is the case when the operator $A$ is bounded, and hence, so is $F(A)$,
\begin{equation*}
\sigma_c(F(A)) = F(\sigma_c(A))
\ \hfill \text{(Spectral Mapping Theorem for Continuous Spectrum)}.
\end{equation*}
\end{enumerate}
\end{cor} 


\begin{cor}[Finer Spectrum Structure]\ \\
Let  $\left\{T(t) \right\}_{t\ge 0}$  be a $C_0$-semigroup $\left\{T(t) \right\}_{t\ge 0}$ (of normal operators) on a complex Hilbert space generated by a normal operator $A$. Then
\begin{align*}
\sigma_p(T(t)) &=e^{t\sigma_p(A)},\ t\ge 0,\ \text{(Spectral Mapping Theorem for Point Spectrum)},\\ 
\sigma_r(T(t))&=\emptyset,\ t\ge 0,\\
\sigma_c(T(t))&=\overline{e^{t\sigma(A)}}\setminus e^{t\sigma_p(A)},\ t\ge 0.
\end{align*}

If furthermore, for some $t>0$, the restriction of the exponential function $e^{t\cdot}$ to $\sigma(A)$ is injective, then
\begin{equation*}
\sigma_c(T(t)) =e^{t\sigma_c(A)}\cup \left(\overline{e^{t\sigma(A)}}\setminus e^{t\sigma(A)}\right).
\end{equation*}

In particular, if the set $e^{t\sigma(A)}$ is closed, as is the case when the operator $A$ is bounded,
\begin{equation*}
\sigma_c(T(t)) =e^{t\sigma_c(A)}
\ \text{(Spectral Mapping Theorem for Continuous Spectrum)}.
\end{equation*}
\end{cor} 

\section{Acknowledgments}

I would like to express my utmost appreciation to Dr.~Michel L. Lapidus for bringing to my attention relevant findings of new monograph \cite{HerLap} while it was still in press  as well as for our interesting and highly stimulating discussions.


\end{document}